\documentclass{amsart}

\usepackage{fancyhdr}
\usepackage{lastpage}
\usepackage{yhmath}

\newcommand{\bdis}{\begin{displaymath}}
\newcommand{\edis}{\end{displaymath}}
\newcommand{\be}{\begin{equation}}
\newcommand{\ee}{\end{equation}}

\newcommand{\mbb}{\mathbb}
\newcommand{\mcal}{\mathcal}

\newcommand{\vp}{\varphi}
\newcommand{\vt}{\vartheta}

\newcommand{\bpsi}{\bar{\psi}}

\newcommand{\zf}{\zeta\left(\frac{1}{2}+it\right)}

\DeclareMathOperator*{\ssum}{\sum\sum}

\DeclareMathOperator{\re}{Re}
\DeclareMathOperator{\im}{Im}

\theoremstyle{definition}

\theoremstyle{remark}
\newtheorem{remark}[]{Remark}

\newtheorem*{mydef11}{{\bf Theorem 1}}

\newtheorem*{mydef12}{{\bf Theorem 2}}

\newtheorem*{mydef13}{{\bf Theorem 3}}

\newtheorem*{mydef14}{{\bf Theorem 4}}

\newtheorem*{mydef21}{{\bf Definition 1}}

\newtheorem*{mydef22}{{\bf Definition 2}}

\newtheorem*{mydef4}{{\bf Corollary}}

\newtheorem*{mydef5A}{{\bf Lemma $\bar{A}$}}
\newtheorem*{mydef5B}{{\bf Lemma $\bar{B}$}}
\newtheorem*{mydef5C}{{\bf Lemma $\bar{C}$}}
\newtheorem*{mydef5alf}{{\bf Lemma $\bar{\alpha}$}}
\newtheorem*{mydef5be}{{\bf Lemma $\bar{\beta}$}}

\numberwithin{equation}{section}

\begin{document}

\title{On Selberg's theorem C in the theory of the Riemann zeta-function}

\author{Jan Moser}

\address{Department of Mathematical Analysis and Numerical Mathematics, Comenius University, Mlynska Dolina M105, 842 48 Bratislava, SLOVAKIA}

\email{jan.mozer@fmph.uniba.sk}


\begin{abstract}
In this paper we obtain new theorems about classes of exceptional sets for the Selberg's theorem C (1942). Our theorems,
as based on discrete method, are not accessible for Karatsuba's theory (1984) since this theory is a continuous theory.
This paper is English version of our paper \cite{8}, the results of our paper \cite{9} are added too.
\end{abstract}

\maketitle

\section{Introduction}

\subsection{}

We use the following notions. Let
\begin{itemize}
 \item[(a)] $\psi(t)$ be a positive increasing to infinity function such that
 \bdis
 \psi(t)\leq \sqrt{\ln t},
 \edis

 \item[(b)] $S$ be the set of values of $t$
 \be \label{1.1}
 t\in [T,T+T^{1/2+\epsilon}]
 \ee
 for which there is at least one zero point of the function
 \be \label{1.2}
 \zf
 \ee
 within the interval
 \be \label{1.3}
 \left(t,t+\frac{\psi(t)}{\ln t}\right),\quad S=S(T,\epsilon,\psi).
 \ee
\end{itemize}

Let us remind that the set of segments (\ref{1.1}) for every small and fixed $\epsilon>0$ is the minimal
set for the Selberg's theory. It is the assertion of the Selberg's C theorem (see \cite{10}, p. 49)
relative to segment (\ref{1.1})
\be \label{1.4}
m(S)\sim T^{1/2+\epsilon},\quad T\to\infty,
\ee
that is, the measure of the set $\bar{S}$ of such values
\bdis
t\in [T,T+T^{1/2+\epsilon}]
\edis
for which there is no zero of the function (\ref{1.2}) in the interval (\ref{1.3}) is
\be \label{1.5}
m(\bar{S})=o(T^{1/2+\epsilon}),\quad T\to\infty.
\ee

\subsection{}

Next, let
\bdis
\{ g_\nu\}
\edis
denote the sequence that is defined by the formula
\bdis
\vt_1(g_\nu)=\frac{\pi}{2}\nu,\quad \nu=1,2,\dots
\edis
(see \cite{6}, $\bar{t}_\nu=g_\nu$, comp. \cite{3}, \cite{4}), where
\bdis
\vt_1(t)=\frac t2\ln\frac{t}{2\pi}-\frac t2-\frac{\pi}{8}.
\edis
\begin{remark}
Since the set
\be \label{1.6}
W=\{ g_\nu:\ g_\nu\in [T,T+T^{1/2+\epsilon}]\}
\ee
is the finite one, then
\be \label{1.7}
m(W)=0.
\ee
Consequently, we have the following: no information is contained in the Selberg's theorem C
(comp. (\ref{1.5}), (\ref{1.7})) about the zeros of odd order of the function (\ref{1.2}) in the
intervals
\bdis
\left( g_\nu,g_\nu+\frac{\psi(g_\nu)}{\ln g_\nu}\right),\ g_\nu\in W
\edis
that is the set $W$ is the exceptional set for the Selberg's theorem C.
\end{remark}

\subsection{}

Now, let us remind the following deep methods of the English mathematicians:
\begin{itemize}
 \item[(a)] continuous method of Hardy-Littlewood (see \cite{1}),
 \item[(b)] discrete method of E.C. Titchmarsh (see \cite{11}).
\end{itemize}

In our papers \cite{6}, \cite{7}, we have constructed a discrete analogue of the
Hardy-Littlewood continuous method (that is, some synthesis of (a) and (b)). Especially, we have
obtained the following estimate (see \cite{7})
\be \label{1.8}
N_0(T+T^{5/12}\psi\ln^3T)-N_0(T)>A(\psi)T^{5/12}\psi\ln^3T,
\ee
where $N_0(T)$ denotes the number of zeros of the function
\bdis
\zf,\quad t\in (0,T],
\edis
and $A(\psi)$ is the constant that depends on choice of $\psi$, for example, if
\bdis
\psi=\ln\ln\ln T,
\edis
then
\bdis
A(\ln\ln\ln T)
\edis
is an absolute constant.

\begin{remark}
We notice explicitly that
\begin{itemize}
 \item[(a)] our improvement of the classical Hardy-Littlewood exponent $\frac 12$
 \bdis
 \frac 12 \longrightarrow \frac{5}{12}
 \edis
 is a $16.6\%$ change after 61 years,

 \item[(b)] the estimate (\ref{1.8}) was the first step on a way to proof of the Selberg's
 hypothesis (see \cite{10}, p. 5, comp. \cite{2}, pp. 37,39).
\end{itemize}
\end{remark}

\begin{remark}
Let us notice that I have sent the manuscripts of my papers \cite{6}, \cite{7} to A.A. Karatsuba
in the beginning of 1981.
\end{remark}

After this analysis, it is clear that our estimate (\ref{1.8}) is in need of a corresponding
analogue of the Selberg's theorem C. Consequently, in this paper we shall prove an analogue
of that theorem for finite set $W_1$ (and also for others) of values
\bdis
g_\nu\in [T,T+T^{5/12}\psi\ln^3T];\quad m(W_1)=0,
\edis
i. e. for the exceptional set in the sense of the Selberg's theorem.

\section{Theorem 1}

\subsection{}

Let (see \cite{6}, (2.5))
\be \label{2.1}
\begin{split}
& \omega=\frac{\pi}{\ln\frac{T}{2\pi}}=\frac{\pi}{2\ln P_0},\
U=T^{5/12}\psi\ln^3T, \\
& \ln T<M<\sqrt[3]{T}\ln T.
\end{split}
\ee
Next, let
\bdis
\bpsi(t)
\edis
be the function of the same kind as $\psi(t)$ and fulfilling the condition
\be \label{2.2}
\frac{\bpsi}{\sqrt[3]{\psi}}=o(1),\ T\to\infty,
\ee
and let
\bdis
G(T,\psi,\bpsi)
\edis
denote  the number of such
\bdis
g_\nu\in [T,T+U]
\edis
that the interval
\be \label{2.3}
(g_\nu,g_\nu+\bpsi(g_\nu))
\ee
contains a zero of the odd order of the function
\bdis
\zf,\ t\in [T,T+U].
\edis
The following theorem holds true.

\begin{mydef11}
\be \label{2.4}
G(T,\psi,\bpsi)\sim \frac{1}{\pi} U\ln T,\ T\to\infty.
\ee
\end{mydef11}

\begin{remark}
Since (see \cite{6}, (8))
\be \label{2.5}
\sum_{T\leq g_\nu\leq T+U}1 \sim \frac{1}{\pi}U\ln T,\ T\to\infty
\ee
then we have by Theorem 1 that for \emph{almost all}
\bdis
g_\nu\in [T,T+U]
\edis
the interval (\ref{2.3}) contains a zero of the odd order of the function
\bdis
\zf.
\edis
\end{remark}

\begin{remark}
Let $N(T)$ denote the number of zeros of the function
\bdis
\zeta(s),\ s=\sigma+it,\ \sigma\in (0,1),\ t\in (0,T].
\edis
It is then true that (comp. \cite{12}, p. 181)
\be \label{2.6}
N(T+U)-N(T)\sim \frac{1}{2\pi}U\ln T,\ T\to\infty.
\ee
Of course, our formula (\ref{2.4}) is not in a contradiction with the formula (\ref{2.6}) since
many of intervals (\ref{2.3}) can intersect.
\end{remark}

\begin{remark}
Let us notice explicitly that also for the theory of Karatsuba giving the estimate
(comp. \cite{2}, p. 39)
\bdis
N_0(T+T^{27/82+\epsilon})-N_0(T)>A(\epsilon)T^{27/82+\epsilon},
\edis
we have that the set of values
\bdis
g_\nu\in [T,T+T^{27/82+\epsilon}]
\edis
is the exceptional set, since the theory is \emph{continuous} as well as the classical
theories of Hardy-Littlewood and Selberg. Consequently, our Theorem 1 is not improvable
also by Karatsuba's theory.
\end{remark}

\subsection{}

Now we give the proof of Theorem 1. The basic point of the proof is the estimate (see \cite{7},
(3.16))
\be \label{2.7}
R<A\frac{U\ln^2T}{M}
\ee
where $R$ denotes the number of such
\bdis
g_\nu^*\in [T,T+U]
\edis
for which the sequence
\bdis
\{ Z(g_\nu^*+k\omega)\}_{k=1}^M
\edis
preserves the sign (comp. \cite{7}, (3.9), (3.11)). Next,
\bdis
\bpsi(g_\nu)\geq \bpsi(T),\ g_\nu\in [T,T+U],
\edis
and
\bdis
\frac{\bpsi(T)}{\omega}\sim \frac 1\pi \bpsi\ln T>\frac{1}{2\pi}\bpsi\ln T\geq
\left[\frac{1}{2\pi}\bpsi\ln T\right]=M_1,
\edis
of course,
\bdis
M_1\in (\ln T,\sqrt[3]{\psi}\ln T)
\edis
(see (\ref{2.1}) - inequalities for $M$ and (\ref{2.2})). Putting $M=M_1$ in (\ref{2.7}) one
obtains
\bdis
R=o(U\ln T).
\edis
Now, the formula (\ref{2.4}) follows from the previous by (\ref{2.5}).

\section{Lemmas about translations $g_\nu\longrightarrow g_\nu(\tau), \tau\in [-\pi,\pi]$}

\subsection{}

Let
\bdis
\{ g_\nu(\tau)\}
\edis
denote the infinite set of sequences which are defined (comp. \cite{5}) by the formula
\bdis
\vt_1[g_\nu(\tau)]=\frac{\pi}{2}\nu+\frac{\tau}{2},\ \nu=1,2,\dots ,\ \tau\in [-\pi,\pi],
\edis
where, of course,
\bdis
g_\nu(0)=g_\nu.
\edis
Now, we shall study how the lemmas from the papers \cite{6}, \cite{7} are sensitive with respect
to the translations
\bdis
g_\nu\longrightarrow g_\nu(\tau), g_\nu\in [T,T+U],\ \tau\in [-\pi,\pi].
\edis
First of all we have (comp. \cite{6}, (22) -- (36)) the following

\begin{mydef5A}
\bdis \begin{split}
& g_{\bar{\nu}_1+p+1}(\tau)=g_{\bar{\nu}_1}(\tau)+\bar{\omega}_0p-\bar{\omega}_0D(p)+
\mcal{O}\left(\frac{U^3}{T^2\ln T}\right), \\
& p=0,1,\dots,N_1-1,
\end{split}
\edis
where (comp. \cite{6}, (11), (12))
\bdis
\begin{split}
 & \bar{\omega}_0=\frac{\pi}{\ln\frac{T}{2\pi}}-\frac{\pi^2}{2}\frac{1}{T\ln^3\frac{T}{2\pi}}-
 \pi\frac{g_{\bar{\nu}_1}(\tau)-T}{T\ln^2\frac{T}{2\pi}}, \\
 & Q=Q(T)=\frac{\pi}{T\ln^2\frac{T}{2\pi}}, \\
 & D(p)=\sum_{q=1}^p \{ 1-(1-Q)^q\}, 1\leq p\leq N_1-1,\ D(0)=0, \\
 & g_{\bar{\nu}_1}(\tau)=\min_{g_\nu(\tau)\in [T,T+U]}\{ g_\nu(\tau)\}, \\
 & g_{\bar{\nu}_1+N_1}(\tau)=\max_{g_\nu(\tau)\in [T,T+U]}\{ g_\nu(\tau)\},\
 \bar{\nu}_1=\bar{\nu}_1(\tau),\ N_1=N_1(\tau),
\end{split}
\edis
and the $\mcal{O}$ is valid uniformly for $\tau\in [-\pi,\pi]$.
\end{mydef5A}

\subsection{}

Next, since (see \cite{6}, (118))
\be \label{3.1}
\bar{\vt}_{1,k}=\vt_1[g_\nu(\tau)+k\omega]=\frac{\pi}{2}\nu+\frac{\tau}{2}+k\omega\ln P_0+
\mcal{O}\left(\frac{MU}{T\ln T}\right),
\ee
then (comp. \cite{6}, (121))
\be \label{3.2}
\begin{split}
 & Z[g_\nu(\tau)+k\omega]\cdot Z[g_\nu(\tau)+l\omega]= \\
 & = 2\ssum_{m,n<P_0}\frac{1}{\sqrt{nm}}
 \cos\{ g_\nu(\tau)\ln\frac nm+k\omega\ln\frac{P_0}{n}-l\omega\frac{P_0}{m}\}+ \\
 & + 2\ssum_{m,n<P_0}\frac{(-1)^\nu}{\sqrt{nm}}
 \cos\{ g_\nu(\tau)\ln(mn)-\tau-k\omega\ln\frac{P_0}{n}-l\omega\frac{P_0}{m}\}+ \\
 & + \mcal{O}\left(\frac{MU}{\sqrt{T}\ln T}\right)+\mcal{O}(T^{-1/12}\ln T),
\end{split}
\ee
and the $\mcal{O}$-estimates in (\ref{3.1}), (\ref{3.2}) are valid uniformly for
$\tau\in[-\pi,\pi]$.

Now, we  put (see (\ref{3.2}), comp. \cite{6}, (16), (17))
\bdis
\begin{split}
& \bar{S}_1(T,U,M,\tau)=\ssum_{m<n<P_0}\frac{1}{\sqrt{mn}}\sum_{T\leq g_\nu(\tau)\leq T+U}\cos\left\{ g_\nu(\tau)\ln\frac nm+\vp_1\right\},
\end{split}
\edis
where
\bdis
\vp_1=k\omega\ln\frac{P_0}{m}-l\omega\ln\frac{P_0}{n},
\edis
and also (comp. \cite{6}, (19), (20))
\bdis
\bar{S}_2(T,U,M,\tau)=\ssum_{m<n<P_0}\frac{1}{\sqrt{mn}}\sum_{T\leq g_\nu(\tau)\leq T+U}
(-1)^\nu\cos\{ g_\nu(\tau)\ln(mn)+\bar{\vp}_2\},
\edis
where
\bdis
\bar{\vp}_2=-k\omega\ln\frac{P_0}{n}-l\omega\ln\frac{P_0}{m}-\tau=\vp_2-\tau.
\edis
Since
\bdis
T\leq g_\nu(\tau)\leq T+U,
\edis
and
\bdis
\begin{split}
 & \bar{S}_2(T,U,M,\tau)=\re\left\{ e^{-i\tau}\ssum_{m,n<P_0}\frac{1}{\sqrt{mn}} \times \right. \\
 & \left. \times \sum_{T\leq g_\nu(\tau)\leq T+U} (-1)^\nu
 \exp\{ i[g_\nu(\tau)\ln(mn)+\vp_2]\}\right\}
\end{split}
\edis
then the following estimates (comp. \cite{6}, (18), (21), (37) -- (93)) hold true

\begin{mydef5B}
\bdis
\bar{S}_1(T,U,M,\tau)=\mcal{O}(MT^{5/12}\ln^3T)
\edis
uniformly for $\tau\in [-\pi,\pi]$.
\end{mydef5B}

\begin{mydef5C}
\bdis
\bar{S}_2(T,U,M,\tau)=\mcal{O}(T^{5/12}\ln^2T)
\edis
uniformly for $\tau\in [-\pi,\pi]$.
\end{mydef5C}

\subsection{}

Let (see \cite{6}, (3))
\bdis
\bar{J}=\bar{J}(T,U,M,\tau)=\sum_{T\leq g_\nu(\tau)\leq T+U}
\left\{\sum_{k=0}^M Z[g_\nu(\tau)+k\omega]\right\}^2.
\edis
Now we obtain by method \cite{6}, (94) -- (127) the following

\begin{mydef5alf}
\bdis
\bar{J}=AMU\ln^2T+o(MU\ln^2T),
\edis
($A>0$ is an absolute constant) uniformly for $\tau\in[-\pi,\pi]$.
\end{mydef5alf}

Next, we have, instead of \cite{7}, (5.1), (5.2), the following
\bdis
\begin{split}
 & 4\cos\bar{\vt}_k\cos\bar{\vt}_l\cos(\bar{\vt}_k-\bar{\vt}_l)= \\
 & = 1+(-1)^{k+l}+(-1)^{\nu+k}\cos\tau+(-1)^{\nu+l}\cos\tau+
 \mcal{O}\left(\frac{MU}{T\ln T}\right), \\
 & -4\cos^2\bar{\vt}_k=-2-2(-1)^{\nu+k}\cos\tau+\mcal{O}\left(\frac{MU}{T\ln T}\right).
\end{split}
\edis
Now, putting (comp. \cite{7}, (3.4), (3.5))
\bdis
\begin{split}
 & \bar{N}=\sum_{T\leq g_\nu(\tau)\leq T+U} |\bar{K}|^2, \\
 & \bar{K}=\sum_{k=0}^M
 \left\{ e^{-i\vt[g_\nu(\tau)+k\omega]}Z[g_\nu(\tau)+k\omega]-1\right\},
\end{split}
\edis
we obtain by method \cite{7}, (4.1) -- (7.3) the following

\begin{mydef5be}
\bdis
\bar{N}=\mcal{O}(MU\ln^2T)
\edis
uniformly for $\tau\in[-\pi,\pi]$.
\end{mydef5be}

\section{Two theorems connected with translations $g_\nu\longrightarrow g_\nu(\tau), \tau\in [-\pi,\pi]$}

\subsection{}

First of all we give (comp. \cite{7}, (2.6) -- (2.8)) the following

\begin{mydef21}
We shall call the segment
\bdis
[g_\nu(\tau)+k(\nu)\omega,g_\nu(\tau)+(k(\nu)+1)\omega]
\edis
where
\bdis
g_\nu(\tau)\in [T,T+U],\ \tau\in [-\pi,\pi],\ 0\leq k(\nu)\leq M_2=[\delta\ln T],\ \delta>1,
\edis
and $k(\nu)\in\mbb{N}_0$ as \emph{the good segment} (comp. \cite{7}, \cite{11}) if
\bdis
Z[g_\nu(\tau)+k(\nu)\omega]\cdot Z[g_\nu(\tau)+(k(\nu)+1)\omega]<0 .
\edis
\end{mydef21}

Next, let
\bdis
G_1(T,U,\delta,\tau)
\edis
denote the number of non-intersecting good segments within the interval $[T,T+U]$. Then we obtain,
similarly to \cite{7}, (3.7), (3.20), the following result

\begin{mydef12}
There are
\bdis
\delta_0>1,\ A(\psi,\delta_0)>0,\ T_0(\psi,\delta_0)>0
\edis
such that
\be \label{4.1}
G_1(T,U,\delta_0,\tau)> A(\psi,\delta_0)U,\ T\geq T_0(\psi,\delta_0)
\ee
for all $\tau\in [-\pi,\pi]$.
\end{mydef12}

\begin{remark}
We notice explicitly that the estimate \cite{7}, (2.9) concerning the number of good segments
(relatively to $\{ g_\nu\}$) is invariant with respect to translations
\bdis
g_\nu\longrightarrow g_\nu(\tau), \ \tau\in [-\pi,\pi],\ g_\nu\in [T,T+U].
\edis
\end{remark}

\subsection{}

Above listed facts make clear that we have obtained a kind of generalization of our
Theorem 1. Namely, let
\bdis
G_2(T,\psi,\bpsi,\tau)
\edis
stand for the number of values
\bdis
g_\nu(\tau)\in [T,T+U]
\edis
such that the interval
\be \label{4.2}
(g_\nu(\tau),g_\nu(\tau)+\bpsi[g_\nu(\tau)])
\ee
contains a zero of the odd order of the function
\bdis
\zf.
\edis
Then the following theorem holds true.

\begin{mydef13}
\bdis
G_2(T,\psi,\bpsi,\tau)\sim \frac 1\pi U\ln T,\ T\to\infty,\ \tau\in[-\pi,\pi].
\edis
\end{mydef13}

\section{Remarks on Selberg's theorems about zeros of function $\zf$}

We have given a discrete commentary to fundamental Selberg's memoir \cite{10} in our paper
\cite{9}. Here we put two results from our paper \cite{9}.

\subsection{}

Let
\be \label{5.1}
H_1\in [a_1,a_2\sqrt{\ln P_0}],
\ee
where
\bdis
a_1=\frac{10}{\pi\epsilon},\quad a_2=a_1\sqrt{\frac 2\pi},\ H_1\in\mbb{N}.
\edis
The origin of (\ref{5.1}) is as follows: we put
\bdis
\omega=\frac{\pi}{2\ln P_0},\ H_1\omega=H,\ \xi=\left(\frac{T}{2\pi}\right)^{\epsilon/10}=
P_0^{\epsilon/5},\ \epsilon\leq \frac{1}{10},
\edis
(see \cite{9}, (9)), and further, we assume that
\bdis
\frac{1}{\ln\xi}\leq H\leq \frac{1}{\sqrt{\ln\xi}},
\edis
(see \cite{9}, (10)).

\begin{mydef22}
We shall call the segment
\bdis
[g_\nu(\tau)+(k(\nu,\tau)-1)\omega,g_\nu(\tau)+k(\nu,\tau)\omega],
\edis
where
\bdis
g_\nu(\tau)\in [T,T+U],\ 1\leq k(\nu,\tau)\leq N_1,
\edis
and $k(\nu,\tau)\in\mbb{N}$ as \emph{the good segment} (see \cite{9}, p. 113) if
\bdis
Z[g_\nu(\tau)+(k(\nu,\tau)-1)\omega]\cdot Z[g_\nu(\tau)+k(\nu,\tau)\omega]<0 .
\edis
\end{mydef22}

Next, let
\bdis
G_3(T,U,H_1,\tau)
\edis
denote the number of non-intersecting good segments within the interval $[T,T+U]$. Then the following
theorem holds true.

\begin{mydef14}
There are
\bdis
\bar{H}_1\in [a_1,a_2\sqrt{\ln P_0}],\ A(\epsilon)>0,\ T_0(\epsilon)>0
\edis
such that
\be \label{5.2}
G_3(T,U,\tau)>A(\epsilon) U\ln T,\ T\geq T_0(\epsilon),
\ee
where, of course,
\bdis
G_3(T,U,\tau)=G_3(T,U,\bar{H}_1,\tau)
\edis
\end{mydef14}

\begin{remark}
Since
\be \label{5.3}
G_3(T,U,\tau)<AU \ln T ,
\ee
then the order of $G_3$ is $U\ln T$ for every fixed $\tau\in [-\pi,\pi]$ (comp. (\ref{5.2}),
(\ref{5.3})). We shall call this property as \emph{generalized Gram's law} for the set of
sequences $\{ g_\nu(\tau)\}$.
\end{remark}

Now, we obtain the following from our Theorem 4.

\begin{mydef4}
\bdis
N_0(T+T^{1/2+\epsilon})-N_0(T)>A(\epsilon)T^{1/2+\epsilon}\ln T,
\edis
i. e. the Selberg's theorem A (see \cite{10}, p. 46, $a=1/2+\epsilon$).
\end{mydef4}

\begin{remark}
Consequently, the generalized Gram's law is the discrete basis of fundamental Selberg's theorem A.
\end{remark}

\begin{remark}
Our Theorem 4 is also not improvable by the Karatsuba theory.
\end{remark}

\begin{remark}
Finally, we mention that our Theorem 3 is valid also for the intervals of the following form
\bdis
\left( g_\nu(\tau),g_\nu(\tau)+\frac{\psi[g_\nu(\tau)]}{\ln g_\nu(\tau)}\right),\
g_\nu(\tau)\in [T,T+T^{1/2+\epsilon}\ln T) ,
\edis
(comp. \cite{9}, p. 112). i.e. for our system of exceptional sets, where
\bdis
a_1\leq \left[\frac{\psi(\tau)}{2\pi}\right]\leq a_2\sqrt{\ln P_0},
\edis
as a discrete analogue of the Selberg's theorem C.
\end{remark}

\thanks{I would like to thank Michal Demetrian for helping me with the electronic version of this work.}

\end{document}